\documentclass[aos,MSNbibl,nameyear,seceqn,dvips]{arximspdf}
\usepackage{mathrsfs}

\doi{10.1214/12-AOS1061} %kopijuoti is PTS
\volume{40}
\issue{6}
\pubyear{2012}
\firstpage{3031}
\lastpage{3049}

\makeatletter
\newcommand{\iint}{\int\!\!\!\int}
\newtheorem{corollary}{Corollary}[section]
\newtheorem{lemma}{Lemma}[section]
\newproclaim{remark}{Remark}[section]
\def\del{\delta}
\def\ve{\varepsilon}
\def\hve{{\hat\varepsilon}}
\def\he{{\hat\varepsilon}}
\def\h{\chi}
\def\vt{\vartheta}
\def\elf{\ell_f}
\def\hvt{{\hat\vt}}
\def\hs{\hat\sigma}
\def\R{\mathbb{R}}
\def\X{\mathcal{X}}
\def\law{\mathfrak{L}}
\def\emp{\mathbb{F}}
\def\hemp{\hat\emp}
\def\1{\mathbf{1}}
\def\argmin{\mathop{\arg\min} }
\def\dj{\delta_j}
\def\di{\delta_i}
\def\d{\delta}
\def\var{\operatorname{Var}}
\def\ndots{\ldots}

\makeatother

\begin{document}
\begin{frontmatter}

\title{The transfer principle: A tool for complete case~analysis}
\runtitle{The transfer principle}

\begin{aug}
\author[A]{\fnms{Hira L.} \snm{Koul}\ead[label=e1]{koul@stt.msu.edu}\ead[label=u1,url]{http://www.stt.msu.edu/\textasciitilde koul}},
\author[B]{\fnms{Ursula U.} \snm{M\"uller}\corref{}\thanksref{t2}\ead[label=e2]{uschi@stat.tamu.edu}\ead[label=u2,url]{http://www.stat.tamu.edu/\textasciitilde uschi}}
\and
\author[C]{\fnms{Anton} \snm{Schick}\thanksref{t3}\ead[label=e3]{anton@math.binghamton.edu}\ead[label=u3,url]{http://www.math.binghamton.edu/anton}}
\thankstext{t2}{Supported by NSF Grant DMS-09-07014.}
\thankstext{t3}{Supported by NSF Grant DMS-09-06551.}
\runauthor{H. L. Koul, U. U. M\"uller and A. Schick}
\affiliation{Michigan State University, Texas A\&M University
and Binghamton~University}
\address[A]{H. L. Koul\\
Department of Statistics and Probability\\
Michigan State University\\
East Lansing, Michigan 48824-1027\\
USA\\
\printead{e1}\\
\printead{u1}}

\address[B]{U. U. M\"uller \\
Department of Statistics \\
Texas A\&M University \\
College Station, Texas 77843-3143\\
USA \\
\printead{e2}\\
\printead{u2}}

\address[C]{A. Schick\\
Department of Mathematical Sciences \\
Binghamton University \\
Binghamton, New York 13902-6000\\
USA\\
\printead{e3}\\
\printead{u3}}
\end{aug}

% HISTORY:
\received{\smonth{8} \syear{2011}}
\revised{\smonth{6} \syear{2012}}

% ABSTRACT
%
\begin{abstract}
This paper gives a general method for deriving limiting distributions of
complete case statistics for missing data models from corresponding
results for the model where all data are observed.
This provides a convenient tool for obtaining
the asymptotic behavior of
complete case versions of established full data methods without
lengthy proofs.

The methodology is illustrated by analyzing three inference procedures for
partially linear regression models with responses missing at random.
We first show that complete case versions of asymptotically efficient
estimators of the slope parameter for the full model are efficient,
thereby solving the problem of
constructing efficient estimators of the slope parameter for this model.
Second, we derive an asymptotically distribution free test for
fitting a normal distribution to the errors.
Finally, we obtain an asymptotically distribution free test for linearity,
that is, for testing that the nonparametric component of these models
is a constant. This test is new both when data are fully
observed and when data are missing at random.
\end{abstract}

% KEYWORDS
% Pirmas kwd is didziosios raides
%
\begin{keyword}[class=AMS]
\kwd[Primary ]{62E20}
\kwd[; secondary ]{62G05}
\kwd{62G10}
\end{keyword}

\begin{keyword}
\kwd{Transfer principle}
\kwd{missing at random}
\kwd{partially linear models}
\kwd{efficient estimation}
\kwd{martingale transform test for normal errors}
\kwd{testing for linearity}
\end{keyword}

\end{frontmatter}

%s1 #&#
\section{Introduction}
\label{sec1}

The basis for regression is a response variable $Y$ and a covariate
vector $X$ which are linked via the formula $Y = r(X) +\ve$,
where $r$ is a regression function and $\ve$ is an error variable.
The analysis is then carried out based on independent
copies $(X_1,Y_1), \ldots, (X_n, Y_n)$ of the pair $(X,Y)$.
We refer to this as the \emph{full model}.
In applications, however, responses may be missing.
The base observation is then a triple $(X, \delta Y, \delta)$, where
$\delta$
is an indicator variable with $E[\delta]=P(\delta=1)>0$.
The interpretation is that for $\delta=1$, one observes the pair $(X,Y)$,
while for $\delta=0,$ one only observes the covariate $X$.
The analysis is now based on independent copies
$(X_1, \delta_1 Y_1, \delta_1),\ndots, (X_n, \delta_n Y_n, \delta_n)$
of the observation $(X, \delta Y, \delta)$.
An accepted way of analyzing such data is by imputing the missing responses.
Here we take a closer look at \textit{complete case analysis}.
This method ignores the %missing
incomplete observations and works
with only the
$N = \sum_{j=1}^n\dj$ completely observed pairs
$(X_{i_1}, Y_{i_1}), \ndots, (X_{i_N}, Y_{i_N})$.
Formally, to each statistic
\[
T_n = t_n(X_1,Y_1, \ndots,
X_n,Y_n)
\]
for the full model there corresponds the \textit{complete case statistic}
\[
T_{c} = t_N(X_{i_1}, Y_{i_1}, \ndots,
X_{i_N}, Y_{i_N}),
\]
which mimics the statistic $T_n$ by treating
$(X_{i_1},Y_{i_1}),\ndots,(X_{i_N},Y_{i_N})$
as if it were a sample of size $N$ from the original setting without
missing data.

Our main result gives a simple and useful method for obtaining the asymptotic
distribution of $T_{c}$.
We show that the limiting distribution of $T_{c}$
coincides with that of $\tilde T_n =
t_n(\tilde X_1,\tilde Y_1, \ndots, \tilde X_n,\tilde Y_n)$
where $(\tilde X_1,\tilde Y_1), \ndots, (\tilde X_n,\tilde Y_n)$
form a random sample drawn from the conditional distribution of
$(X,Y)$, given
$\delta=1$; see Remark~\ref{rem4}.
This can be used as follows. One typically knows the limiting distribution
$\law(Q)$ of $T_n$ under each joint distribution $Q$ of $X$ and $Y$
belonging to some model. If the distribution $\tilde Q$ of $(\tilde
X,\tilde Y)$
belongs to this model, then the limiting distribution of
the complete case statistic is $\law(\tilde Q)$.
We refer to this as the \textit{transfer principle}.
It provides a convenient tool for obtaining %to obtain
the asymptotic behavior of
complete case versions of established full data methods without (reproducing)
lengthy proofs.

Of special interest are statistics $T_n$ that are \emph{asymptotically linear}
for a functional~$T$ from a class $\mathcal{Q}$ of distributions into
$\R^m$
in the sense that if $X$ and $Y$ have joint distribution $Q$ belonging
to the model
$\mathcal{Q}$, then the expansion
\[
T_n = T(Q)+ \frac1n \sum_{j=1}^n
\psi_Q(X_j,Y_j)+ o_P
\bigl(n^{-1/2}\bigr)
\]
holds. Here $\psi_Q$ is a measurable function into $\R^m$
such that $E[\psi_Q(X,Y)]=0$ and $E[\|\psi_{Q}(X,Y)\|^2]$
is finite when $X$ and $Y$ have joint distribution $Q$.
Here and below $\|\cdot\|$ denotes the Euclidean norm.
The function $\psi_Q$ is commonly called an \emph{influence function}.
From the above expansion we obtain that $n^{1/2}(T_n-T(Q))$ is asymptotically
normal with the zero vector as mean and with dispersion matrix
$\Sigma(Q)=E[\psi_Q(X,Y)\psi_Q^{\top}(X,Y)]$.
If $\tilde Q$ belongs to the model $\mathcal{Q}$,
then we have the expansion
\[
\tilde T_n = T(\tilde Q) + \frac1n \sum
_{j=1}^n \psi_{\tilde
Q}(\tilde
X_j,\tilde Y_j)+ o_P\bigl(n^{-1/2}
\bigr),
\]
and obtain from our main result
that
\[
T_{c} = T(\tilde Q)+ \frac{1}{N} \sum
_{j=1}^n\delta_j \psi_{\tilde Q}(X_j,Y_j)+
o_P\bigl(n^{-1/2}\bigr),
\]
see Remark~\ref{rem5}.
From this we immediately derive the expansion
\[
T_{c}= T(\tilde Q)+ \frac{1}{nE[\delta]} \sum
_{j=1}^n\delta_j \psi_{\tilde
Q}(X_j,Y_j)
+ o_P\bigl(n^{-1/2}\bigr).
\]
Thus, if $\tilde Q$ belongs to the model $\mathcal Q$
and $T(\tilde Q)$ equals $T(Q)$, then $T_{c}$ is asymptotically linear
in the model with missing data with influence function $\tilde\psi$,
where
\[
\tilde\psi(X,\delta Y, \delta)= \frac{\delta}{E[\delta]} \psi_{\tilde Q}(X,Y).
\]
We refer to this as the
\emph{transfer principle for asymptotically linear statistics}.
It yields that $n^{1/2}(T_{c}-T(Q))$ is asymptotically normal
with the zero vector as mean and with dispersion
matrix $(1/E[\delta])\Sigma(\tilde Q)$.

The key to a successful application of the transfer principle is the
condition $T(\tilde Q)=T(Q)$. Under this condition, $n^{1/2}$-consistency
carries over to the complete case statistic. If this condition is not met,
the complete case statistic will be biased for estimating $T(Q)$.
% and should not be used.

For our illustration of the transfer principle we consider the
important case
where the response $Y$ is \emph{missing at random} (MAR).
This means that the indicator~$\delta$ is conditionally independent
of $Y$, given $X$, that is,
\[
P(\delta=1|X,Y) = P(\delta=1|X) = \pi(X)\qquad \mbox{a.s.}
\]
This is a common assumption and reasonable in many applications
[see \citet{lit02}, Chapter 1].
This model is referred to as the \emph{MAR model}.

It is well known that the complete case analysis does not always
perform well and that an approach which imputes missing values often
has better statistical properties.
See, for example, Chapter 3 of \citet{lit02} for examples
where using the complete case approach results in bias or a loss of
precision, due to the loss of information.
For a discussion of various imputing methods we again refer to
\citet{lit02}, and also to M\"uller, Schick and Wefelmeyer (\citeyear{msw06}),
who propose efficient estimators for various regression settings which
impute missing and non-missing responses.

Although complete case analysis can lead to the above-mentioned problems,
there are situations where it provides useful and optimal inference procedures.
\citet{efr11}, for example, considers nonparametric
regression with responses missing at random. He shows that his complete case
estimator of the regression function is optimal in the sense that it satisfies
an asymptotic sharp minimax property. M\"uller (\citeyear{mue09})
demonstrates efficiency of a complete case estimator for the parameter vector
in the nonlinear regression model.

For simplicity and clarity, we illustrate the above transfer principle
using the partially linear regression model.
In this model
the response $Y$ is linked to covariates $U$ and $V$ via the relation
%
%e1.1 #&#
\begin{eqnarray}
\label{rho} Y = \vt^\top U + \rho(V) + \ve,
\end{eqnarray}
with $\vt$ an unknown $m$-dimensional vector and $\rho$ an unknown
twice continuously differentiable function.
The error $\ve$ is assumed to have mean zero, finite variance~$\sigma^2$
and a density $f$, and is independent
of the covariates $(U,V)$, where the random vector $U$ has dimension $m$
and the random variable $V$ takes values in the compact interval $[0,1]$.
Throughout this paper, we impose the following conditions on the joint
distribution $G$ of $U$ and $V$:

\begin{longlist}[(G1)]
\item[(G1)]
The covariate $V$ has a density that is bounded and bounded away from
zero on $[0,1]$.
\item[(G2)]
The covariate vector $U$ satisfies $E[\|U\|^2] < \infty$ and the matrix
\[
W_G=E\bigl[\bigl(U-\mu_G(V)\bigr) \bigl(U-
\mu_G(V)\bigr)^{\top}\bigr]
\]
is positive definite, where $\mu_G(V)=E[U|V]$.
\end{longlist}

The requirement involving $W_G$ is needed to identify the parameter
$\vt$.

One important problem is the \textit{efficient} estimation of the regression
parameter~$\vt$ in (\ref{rho}). This is addressed in our first
illustration of the transfer principle below.
The crucial condition for a successful application of the transfer principle,
$T(\tilde Q) = T(Q)$, is satisfied in this case and, more generally,
also for functionals of the triple $(\vt, \rho, f)$.
The MAR assumption and the independence of $\ve$ and $(U,V)$
imply that $\ve$ and $(U,V,\delta)$ are independent.
Hence, the regression parameters $\vt$ and $\rho$ and the error
density $f$
stay the same when conditioning on $\delta=1$. Only the covariate distribution
$G$ changes to $\tilde G$, the conditional distribution of
$(U,V)$ given $\delta=1$.
This argument suggests that inference about the triple $(\vt, \rho, f)$
should be carried out using a complete case analysis, because
the complete case observations are \textit{sufficient} for
$(\vt, \rho, f, \tilde G)$ since they carry all the information
about these parameters.
The covariate pair $(U,V)$ alone, on the other hand,
has no information on $(\vt, \rho, f)$, and hence has no bearing on the
inference about these parameters when the response $Y$ is missing at random.
The same reasoning also applies to general semiparametric regression
models: inference about the regression function and the error distribution
should be based on the complete cases only.

In order to obtain an efficient estimator for $\vt$ we must
assume that the error density $f$ has finite Fisher information for location.
This means that $f$ is absolutely continuous with a.e. derivative $f'$
such that $J_f=\int\elf^2(x)f(x) \,dx$ is finite, where $\elf=-f'/f$
is the score function for location.
Efficient estimators of $\vt$ in the full model are characterized by
the stochastic expansion
\[
\hvt_n = \vt+ \frac1n \sum_{j=1}^n
(J_f W_G)^{-1}\bigl(U_j-
\mu_G(V_j)\bigr)\elf(\ve_j) +
o_P\bigl(n^{-1/2}\bigr);
\]
see, for example, \citet{sch93}. Because of the structure of the
MAR model
%mentioned
introduced
above, the transfer principle for asymptotically linear statistics
yields that the complete case version $\hvt_{c}$ of an
efficient estimator satisfies the expansion
%
%e1.2 #&#
\begin{equation}
\label{eff-mar}\qquad \hvt_{c} = \vt+ \frac1n \sum
_{j=1}^n \frac{\delta_j}{E[\delta]} (J_f
W_{\tilde G})^{-1} \bigl(U_j- \mu_{\tilde G}(V_j)
\bigr) \elf(\ve_j) + o_P\bigl(n^{-1/2}\bigr).
\end{equation}
This of course requires that $\tilde G$ satisfies the properties (G1)
and (G2). This is the case when $\pi$ is bounded away from zero;
see Remark~\ref{rem6a}. Here $\pi(X) = \pi(U,V) = P(\delta= 1|U,V)$.

Although several estimators exist which are efficient in the full
partially linear model, to our knowledge no efficient estimators
have so far been constructed for the corresponding MAR model.
We show in Section~\ref{sec3} that the expansion (\ref{eff-mar})
of $\hvt_c$ characterizes asymptotically
\textit{efficient} estimators of $\vt$
in the MAR model. This means that complete
case versions of efficient estimators in the full model remain efficient
in the MAR model (under appropriate conditions). This result in turn solves
the important problem of constructing efficient estimators for
$\vartheta$
in the partially linear MAR model.
For constructions of efficient estimators in the full model~(\ref{rho}),
we refer the reader to \citet{Cuzick92}, Schick (\citeyear
{sch93}, \citeyear{Sch96a}),
\citet{BZ97} and Forrester et al. (\citeyear{for03}).
Some of these constructions require smoothness assumptions on $\mu_G$.
Then the validity of (\ref{eff-mar}) requires the same smoothness assumptions
on $\mu_{\tilde G}$.

The above method of constructing efficient estimators for the
finite-dimen\-sional parameter also yields efficient estimators
in other semiparametric regression MAR models. The influence function of
the complete case version of an estimator efficient for the full model
is given by the transfer principle for asymptotically linear
estimators. One then only needs to show that this influence function
is the efficient influence function for the MAR model. The latter
can be done by mimicking the results in Section~\ref{sec3}.
There we sketch this approach for the partially linear model
with additive $\rho$ and for a single index model.
M\"uller (\citeyear{mue09}) has calculated the efficient influence function
for the regression parameter in a nonlinear regression model.
Using the transfer principle, one sees that the efficient influence function
equals the influence function of the complete case version
of an efficient estimator for the full model.
This provides a simple derivation of efficient estimators in her model.

We believe that the above \emph{efficiency transfer} is valid for
the estimation of other characteristics in the MAR model (\ref{rho}).
We expect that the efficiency transfer generalizes to the estimation of
(smooth) functionals of the triple $(\vt, \rho, f)$.
This includes as important
special cases the estimation of the \mbox{error} distribution function,
the error variance and other characteristics of $f$ such as quantiles
and moments of~$f$. However, further research is needed to crystallize
the issues involved.

Next we illustrate the transfer principle on %some
goodness-of-fit and
lack-of-fit tests. There is a vast literature on goodness-of-fit tests
for fitting an error distribution and lack-of-fit tests for fitting
a regression function in fully observable regression models.
See, for example, \citet{hart97} and the review article by
\citet{kou06},
and the references
therein. Here we shall discuss two important examples for the
MAR regression models. One pertains to fitting a parametric
distribution to the error distribution in (\ref{rho}) and the other
to testing whether $\rho$ in the model (\ref{rho}) is a constant or not.
In both examples the proposed tests are complete case analogs of full
model tests that are \textit{asymptotically distribution free},
that is, the limiting distribution of the test statistic under the null
hypothesis is the same for all members of the null model being fitted.
Due to the transfer principle, the same conclusion continues to hold
for the proposed tests for the MAR model (\ref{rho}).

First, consider the goodness-of-fit testing problem in the model (\ref{rho})
and the null hypothesis $H_0\dvtx\ve\sim N(0, \sigma^2)$,
for some unknown $0 < \sigma^2 < \infty$.
For the full model a residual-based test of this hypothesis was
introduced by M\"uller, Schick and Wefelmeyer (\citeyear{msw11}) (MSW) adapting
a martingale transform test of \citet{khkou09} for fitting
a parametric family of error distributions in nonparametric regression.
In (\ref{rho}), the residuals
are of the form $\hve_j = Y_j - \hvt^\top U_j - \hat\rho(V_j)$,
where $\hvt$ is a $\sqrt{n}$-consistent estimator of $\vt$
and $\hat\rho$ is a nonparametric estimator of $\rho$, such as
a local smoother based on the covariates $V_j$ and the modified
responses $Y_j-\hvt^\top U_j$, or a series estimator.
Let $\hat\sigma= (\sum_{j=1}^n\he_j^2/n )^{1/2}$ denote the
estimator
of the standard deviation $\sigma$ and $\hat Z_j= \he_j/\hat\sigma$,
$j=1,\ndots,n,$ denote the standardized residuals.
The test statistic of MSW is then
\[
T_n=\sup_{t\in\R}\Biggl | \frac{1}{\sqrt{n}} \sum
_{j=1}^n \bigl\{\1 [\hat Z_j\le t] - H(
\hat Z_j\wedge t)h(\hat Z_j) \bigr\} \Biggr|
\]
for some known functions $h$ and $H$ related to the standard normal
distribution function and its derivatives; see Section~\ref{sec4},
equation (\ref{hH}).
Here we work with a series estimator of $\rho$, which is
discussed in Section 4 of MSW. This requires no additional assumptions.
The test based on $T_n$ is asymptotically distribution free,
because under the null hypothesis, $T_n$ converges in distribution to
%
%e1.3 #&#
\begin{equation}
\label{z} \zeta= \sup_{0\le t\le1} \bigl|B(t)\bigr|,
\end{equation}
where $B$ is a standard Brownian motion.
Due to the transfer principle, the complete case version $T_{c}$ of the
above $T_n$ has the same limiting distribution under the null hypothesis.
Hence, the null hypothesis
is rejected if $T_c$ exceeds the upper $\alpha$ quantile of
the distribution of $\zeta$.
See Section~\ref{sec4}, equation (\ref{tc}), and the discussion around
it for a detailed description of the complete case variant $T_{c}$ of
the above $T_n$. From the discussion on optimality of this test in
\citet{khkou09} and the transfer principle, it follows that the
test based
on $T_c$ will generally be more powerful than the complete case test
based on the Kolmogorov--Smirnov statistic.

Finally, we consider testing whether $\rho$ is constant
within the partially linear model, that is, we
suppose that the partially linear model (\ref{rho})
holds true and test whether the regression function is in fact linear.
%Next, consider testing whether $\rho$ is constant in (\ref{rho}).
Here we adapt an approach by \citet{stu08}
for testing a general parametric model in nonparametric regression,
which is based on a weighted residual-based empirical process.
For the full model this suggests a test statistic of the form
\[
\sup_{t\in\R} \Biggl| \frac{1}{\sqrt{n}} \sum_{j=1}^n
W_j 1[\hve_{j0}\le t] \Biggr|,
\]
where $\hve_{j0}$ are the residuals under the
null hypothesis obtained by regressing the responses $Y_j$ on
the covariates $U_j$ including an intercept,
and where $W_j$ are normalized versions of the residuals obtained from
regressing $\h(V_j)$ on the covariates $U_j$ including an intercept,
for a suitably chosen function $\h$.
The asymptotic null distribution of this test is that of
\[
\zeta_0= \sup_{0\le t \le1} \bigl|B_0(t)\bigr|,
\]
where $B_0$ denotes a standard Brownian bridge.
This is the first test for this problem in the case of fully
observed data. The transfer principle immediately shows that the complete
case variant of this test described at (\ref{tc2}) has the same
limiting distribution.

The literature on lack-of-fit testing in the regression model
when responses are missing at random is scant.
\citet{sw09} establish asymptotic distributional properties
of some tests based on marked residual empirical processes for
fitting a parametric model to the regression function when data are
imputed using the inverse probability method.
\citet{swd09} derive tests to check the hypothesis that
the partially linear model (\ref{rho}) is appropriate,
based on data which are ``completed''
by imputing estimators for the responses.
These tests are compared with tests that ignore the missing data
pairs. Gonz\'alez-Manteiga and
P\'erez-Gonz\'alez (\citeyear{gp06}) use imputation to complete the data.
They derive tests about linearity of the regression function
in a general nonparametric regression model.
Their test is similar to the above test for the last example.
The last two papers report simulation results
that support the superiority of these methods over a selected
complete case method. However, one can verify that
the first test statistic in Sun, Wang and Dai (\citeyear{swd09}) is
asymptotically equivalent to a complete case
statistic in their case 3, and this complete case statistic
should thus result in an equivalent test.
Finally, \citet{li11} uses imputation together with the minimum
distance methodology of \citet{kouni04} to propose
tests for fitting a class of parametric models to the
regression function that includes polynomials.

This article is organized as follows.
Section~\ref{sec2} provides the theory for the transfer principle.
The key is Lemma~\ref{lm1}, which calculates the
explicit form of the distribution of a complete case statistic.
In Section~\ref{sec3} we show that the influence function of the
complete case version of an efficient estimator of $\vt$
in the full data partially linear model is the efficient influence function
for estimating $\vt$ in the MAR model.
Similar results are sketched for a partially linear additive model
(see Remark~\ref{rem6a}) and a single index model (see Remark~\ref{rem6b}).
Section~\ref{sec4} discusses the test for normality of the errors
for the MAR model, and derives expansions for the complete case
residual-based empirical distribution function.
In Section~\ref{sec5} we provide details for the complete case version
of the second test about the nonparametric part $\rho$ in~(\ref{rho})
being constant.

%s2 #&#
\section{Distribution theory for general complete case statistics}
\label{sec2}

In this section we derive the exact distribution of a complete case statistic
in a general setting. Let $(\X,\mathscr{A})$ be a measurable space,
and, for each integer $k$, let $t_k$ be a measurable function from $\X^k$ into
$\R^m$. Let $(\d_1, \xi_1), (\d_2, \xi_2),\ndots$ be independent copies
of $(\d, \xi)$, where $\d$ is Bernoulli with
parameter $p>0$ and $\xi$ is a $\X$-valued random variable.
We denote the conditional distribution of $\xi$, given $\d= 1$ by
$\tilde Q$.
Let $\tilde\xi_1,\tilde\xi_2,\ndots$ be independent $\X$-valued random
variables with common distribution $\tilde Q$.
Denote the distribution of $t_n(\tilde\xi_1,\ndots,\tilde\xi_n)$
by $R_n$.
Then, for any Borel set $B$,
\begin{eqnarray*} R_n(B)&=& \tilde Q^n(t_n
\in B)=P\bigl(t_n(\tilde\xi_1,\ndots,\tilde
\xi_n)\in B\bigr)
\\
&=&P\bigl(t_n(\xi_1,\ndots,\xi_n)\in B |
\delta_1=1,\ndots,\delta_n=1\bigr).
\end{eqnarray*}

By a \emph{complete case statistic associated with the sequence $(t_n)$}
we mean a statistic $T_{c,n}$ of the form
\[
T_{c,n} = \sum_{A \subset\{1, \ndots,n \}} \biggl\{ \prod
_{i \in A} \d_i \biggr\} \biggl\{ \prod
_{i \notin A} (1-\d_i) \biggr\} t_{|A|}
\bigl(\xi^A\bigr),
\]
where $t_0(\xi^\varnothing)$ is a constant, $|A|$ denotes the
cardinality of $A$
and $\xi^A$ is the vector $(\xi_{i_1}, \ndots, \xi_{i_k})$
with $i_1 < \cdots< i_k$ the elements of the non-empty subset
$A \subset\{1, \ndots, n\}$.
Note that the product
$[\prod_{i \in A}\d_i ][\prod_{i \notin A} (1-\d_i)]$ is the
indicator function of the event
$\{\di=1, i\in A\}\cap\{\di=0, i\notin A\}$.
Thus, $T_{c,n}$ equals $t_{|A|}(\xi^A)$ on this event.
It is now clear that $T_{c,n}$ depends on the indicators $\d_1,\ndots
,\d_n$
and only those observations $\xi_i$ for which $\di=1$.

%re2.1 #&#
\begin{remark}\label{rem1}
For a measurable function $\psi$ defined on $\X$, we define the sequence
$(\bar\psi_n)$ by $\bar\psi_n(x_1,\ndots,x_n)=(\psi(x_1)+\cdots
+\psi(x_n))/n$.
The complete case statistic associated with $(\bar\psi_n)$ is
$\sum_{j=1}^n \delta_j\psi(\xi_j)/\sum_{j=1}^n \delta_j$.
\end{remark}

%re2.2 #&#
\begin{remark}\label{rem2}
If $T_{c,n}$ is a complete case statistic associated with $(t_n)$ and
$\alpha$ is a real number, then $(\sum_{j=1}^n \delta_j)^{\alpha}
T_{c,n}$ is a complete case
statistic associated with the sequence $(n^{\alpha}t_n)$.
\end{remark}

For the remainder of this section $T_{c,n}$ denotes a complete case
statistic associated with $(t_n)$ and $H_n$ its distribution.
The next lemma calculates $H_n$ explicitly.

%le2.1 #&#
\begin{lemma}\label{lm1}
For every Borel subset $B$ of $\R^m$, we have
\[
H_n(B)=P(T_{c,n} \in B) = \sum
_{k=0}^n \pmatrix{n \cr k} p^k
(1-p)^{n-k} R_k(B),
\]
with $R_0(B)=1[t_0(\xi^{\varnothing})\in B]$.
\end{lemma}

\begin{pf}
Conditioning on $\d_1,\ndots,\d_n$ yields the identity
\[
P(T_{c,n}\in B)= E\bigl[P(T_{c,n} \in B |\d_1,
\ndots, \d_n)\bigr]
\]
and, thus,
\[
H_n(B)=\sum_{A \subset\{1,\ndots,n\}} p^{|A|}
(1-p)^{n-|A|} H(A,B),
\]
where
\begin{eqnarray*}H(A,B) &=&P(T_{c,n} \in B |
\d_i=1, i \in A, \d_j =0, j \notin A)
\\
&=& P\bigl(t_{|A|}\bigl(\xi^{A}\bigr) \in B |
\d_i=1, i \in A, \d_j =0, j \not \in A\bigr)
\\
&=& \tilde Q^{|A|} (t_{|A|} \in B) = R_{|A|}(B)
\end{eqnarray*}
for non-empty $A$, while $H(\varnothing,B)= R_0(B)$.
The desired result is now immediate.
\end{pf}

%re2.3 #&#
\begin{remark}\label{rem3} Lemma~\ref{lm1} has the following interpretation.
The statistic $T_{c,n}$ has the same distribution as
$t_K(\tilde\xi_1,\ndots,\tilde\xi_K)$,
where $K$ is a binomial random variable with parameters $n$ and $p$,
independent of $\tilde\xi_1,\tilde\xi_2,\ndots.$\vadjust{\goodbreak}
\end{remark}

From the lemma we immediately obtain the following results.

%co2.1 #&#
\begin{corollary} %\label{c1}
The following statements hold:
\begin{longlist}[(a)]
\item[(a)]
If the sequence $(R_n)$ is tight, so is the sequence $(H_n)$.
\item[(b)]
If $R_n$ converges weakly to some limit $L$, then $H_n$ converges weakly
to the same limit $L$.
\item[(c)]
If $R_n$ converges weakly to point mass at $0$, then $T_{c,n}$
converges in
probability to zero.
\end{longlist}
\end{corollary}

%re2.4 #&#
\begin{remark}\label{rem4}
Recall that $R_n$ is the distribution of
$t_n(\tilde\xi_1,\ndots,\tilde\xi_n)$.
Thus, by (b), the limiting distribution of $T_{c,n}$ equals the limiting
distribution of $t_n(\tilde\xi_1,\ndots, \tilde\xi_n)$.
This provides the basis for the transfer principle.
\end{remark}

%re2.5 #&#
\begin{remark}\label{rem5}
Let $\psi$ and $\bar\psi_n$ be as in Remark~\ref{rem1}.
Set $N= \sum_{j=1}^n \delta_j$. Then
\[
S_{c,n} = \sqrt{N} \Biggl(T_{c,n}-\frac{1}{N} \sum
_{j=1}^n \delta_j \psi(
\xi_j) \Biggr)
\]
is a complete case statistic associated with
$s_n = (n^{1/2}(t_n - \bar\psi_n))$.
Suppose that
\[
s_n(\tilde\xi_1,\ndots,\tilde\xi_n) =
n^{1/2} \Biggl(t_n(\tilde\xi_1,\ndots,\tilde
\xi_n)- \frac1n \sum_{j=1}^n
\psi (\tilde\xi_j) \Biggr) = o_P(1).
\]
Then, by (c), we have
\[
S_{c,n} = \sqrt{N} \Biggl(T_{c,n}-\frac{1}{N} \sum
_{j=1}^n \delta_j \psi(
\xi_j) \Biggr) = o_P(1)
\]
and, consequently,
\[
T_{c,n} = \frac{1}{N} \sum_{j=1}^n
\delta_j \psi(\xi_j) + o_P
\bigl(n^{-1/2}\bigr).
\]
This is the basis for the transfer principle for asymptotically
linear statistics.
\end{remark}

%s3 #&#
\section{Efficiency considerations for the partially linear MAR model}
\label{sec3}

Here we shall show that the expansion (\ref{eff-mar}) characterizes
efficient estimators in the partially linear MAR model.
For this we only need to show that the influence function
appearing in (\ref{eff-mar}) is the efficient influence function
for estimating $\vt$ in this model.
We formulate this as the main result of this section; see Lemma~\ref{lm2}.
By the discussion in the \hyperref[sec1]{Introduction}, we
must require that the conditional distribution $\tilde G$
of $(U,V)$, given $\delta=1$, satisfies the assumptions (G1) and (G2).
This is crucial for the transfer principle to apply, and
holds if the function $\pi$ is bounded away from zero, as we shall show
first.

%re3.1 #&#
\begin{remark}\label{rem6a}
Consider the conditional distribution $\tilde G$ of $(U,V)$
given $\delta=1$. Then $\tilde G$ satisfies the properties (G1)
and (G2) if $\pi$ is bounded away from zero:
it is easy to check that $\tilde G$ has density $\tilde\pi$
with respect to $G$, where $\tilde\pi(U,V)=\pi(U,V)/E[\delta]$.
If $\tilde\pi\ge\eta$ for some positive constant $\eta$, then
\[
\eta\int|h| \,dG \le\int|h| \,d\tilde G \le\int|h| \,dG /E[\delta]
\]
for all $h\in L_1(G)$ and, therefore,
\begin{eqnarray*}
a^{\top}W_{\tilde G}a &=&
\int\bigl|a^{\top} \bigl(u-\mu_{\tilde G}(v)\bigr)\bigr|^2 \,d\tilde
G(u,v) \ge\eta\int\bigl|a^{\top}\bigl(u-\mu_{\tilde G}(v)
\bigr)\bigr|^2 \,dG(u,v)
\\
&\ge&\eta\int\bigl|a^{\top} \bigl(u-\mu_G(v)\bigr)\bigr|^2
\,dG(u,v) =\eta a^{\top}W_G a \qquad\mbox{for all } a \in
\R^m.
\end{eqnarray*}
From these inequalities we conclude that
$\tilde G$ inherits the properties (G1) and (G2) from $G$ if
$\pi$ is bounded away from zero.
\end{remark}

%le3.1 #&#
\begin{lemma}\label{lm2}
Suppose the model (\ref{rho}) holds with $\rho$ being twice continuously
differentiable and error density having finite Fisher information for location.
Also assume $\pi$ is bounded away from zero.
%Suppose the assumptions for model (\ref{rho}) from the Introduction
%are satisfied and $\pi$ is bounded away from zero.
Then an efficient estimator of the parameter~$\vt$
in the MAR model is characterized by (\ref{eff-mar}). As a consequence,
the complete case version of an efficient estimator of the parameter~$\vt$
in the full model is efficient for the MAR model.
\end{lemma}

\begin{pf}
We rely heavily on the calculations in M{\"u}ller, Schick and
Wefelmeyer (\citeyear{msw06}).
The authors considered the general missing data problem
with base observation $(X,\delta Y, \delta)$ where $X$ and $Y$
do not have to follow a regression model.
They expressed the joint distribution $P$ of $(X,\delta Y, \delta)$ via
\[
P(dx,dy,dz)= G(dx) B_{\pi(x)}(dz) \bigl(z Q(x,dy) +(1-z)
\Delta_0(dy)\bigr)
\]
in terms of the distribution $G$ of $X$,
the conditional probability $\pi(x)$ of $\delta=1$ given $X=x$,
and the conditional distribution $Q(x,dy)$ of $Y$ given $X=x$.
Here $B_p$ denotes the Bernoulli distribution with parameter $p$
and $\Delta_t$ the Dirac measure at $t$.
They showed that the tangent space is the sum of the orthogonal spaces
\begin{eqnarray*}
 T_1 &=&\bigl\{ u(X)\dvtx u\in\mathscr{U}
\bigr\},\qquad  T_2 = \bigl\{ \delta v(X,Y)\dvtx v\in\mathscr{V}\bigr\},
\\
T_3 &=&\bigl\{ \bigl(\delta-\pi(X)\bigr)w(X)\dvtx w\in\mathscr{W}
\bigr\}.
\end{eqnarray*}
Here, the set $\mathscr{U}$ consists of all real-valued functions $u$
satisfying $\int u \,dG=0$, $\int u^2 \,dG<\infty$ and for which there is
a sequence $G_{nu}$ of distributions fulfilling the model assumptions
on $G$ and
\[
\int \biggl(n^{1/2} \bigl(dG_{nu}^{1/2}-dG^{1/2}
\bigr)-\frac{1}{2} u\, dG^{1/2} \biggr)^2 \to0.
\]
The set $\mathscr{W}$ consists of real-valued functions $w$ with
the property $\int w^2\pi(1-\pi) \,dG<\infty$ for which there is
a sequence $\pi_{nw}$ satisfying the model assumptions on $\pi$
such that
\[
\int \biggl(n^{1/2} \bigl(dB_{\pi_{nw}(x)}^{1/2}-dB_{\pi(x)}^{1/2}
\bigr) - \frac{1}{2}\bigl(\cdot-\pi(x)\bigr) \,dB^{1/2}_{\pi(x)}
\biggr)^2 G(dx) \to0.
\]
Finally, the set $\mathscr{V}$ consists of functions $v$ with the
properties $\int v(x,y)Q(x,dy)$ $=0$ for all $x$
and $\int v^2(x,y)G(dx)Q(x,dy)<\infty$,
and for which there is a sequence $Q_{nv}$ satisfying the model
assumptions on $Q$ and
\[
\iint \biggl(n^{1/2} \bigl(dQ_{nv}^{1/2}(x,
\cdot)-dQ^{1/2}(x,\cdot)\bigr) - \frac12 v(x,\cdot) \,dQ^{1/2}(x,
\cdot) \biggr)^2 G(dx) \to0.
\]

In the partially linear regression model (\ref{rho})
we have $X=(U,V)$ and
\[
Q(x,dy)= %Q_{\vt,\rho,f}(x,dy)= f(y-\vt^{\top}U-\rho(V)) \,dy,
Q_{\vt,\rho,f}(u,v,dy)= f\bigl(y-\vt^{\top}u-
\rho(v)\bigr) \,dy,
\]
where the density $f$ has finite Fisher information for location,
$\vt$ belongs to $\R^m$ and $\rho$ is a smooth function.
For this model $\mathscr{V}$ consists of the functions
\[
a^{\top}U\elf(\ve)+ b(V)\elf(\ve) +c(\ve)
\]
with $a\in\R^m$, $E[b^2(V)]<\infty$
and $c\in L_2(F)$ with $\int c(y) \,dF(y)=0$
and\break $\int c(y)y \,dF(y)=0$.
Since we are interested in estimating the finite-dimensional
parameter $\vt$, we introduce the functional
\[
\kappa(G,Q_{\vt,\rho,f},\pi)=\vt.
\]
Now consider
\[
g(X,\delta Y, \delta)= \delta\bigl(U-\mu_1(V)\bigr) \elf(\ve)
\]
with $\mu_1(V)=E(U|V,\delta=1)$. Then the coordinates of
$g(X,\delta Y, \delta)$ belong to $\mathscr{V}$.
Thus, we have $E[g(X,\delta Y, \delta) u(X)]=0$ and
$E[g(X,\del Y, \delta)(\delta-\pi(X))w(X)]$ $=0$.
Note that $\ve$ and $(\delta,X)$ are independent,
and that we have $E[\elf(\ve)]=0$ and $E[\elf^2(\ve)]=J_f$.
Using this and the definition of $\mu_1$, we calculate
\begin{eqnarray*}
&&E\bigl[g(X, \del Y, \del)\delta
\bigl(a^{\top}U\elf(\ve)+b(V)\elf(\ve )+c(\ve)\bigr)\bigr]
\\
&&\qquad= E\bigl[\delta\bigl(U-\mu_1(V)\bigr) \bigl(U^{\top}a +
b(V)\bigr)\bigr] J_f + E\bigl[\delta\bigl(U-\mu_1(V)
\bigr)\bigr] E\bigl[\elf(\ve)c(\ve)\bigr]
\\
&&\qquad=E\bigl[\delta\bigl(U-\mu_1(V) \bigr) \bigl(U-
\mu_1(V)\bigr)^\top\bigr] a J_f.
\end{eqnarray*}
From this we can conclude that the functional $\kappa$
is differentiable with canonical gradient $g_*(X, \del Y, \del)$
of the form
\begin{eqnarray*}
&& \delta \bigl( J_f E\bigl[\delta
\bigl(U-\mu_1(V)\bigr) \bigl(U-\mu_1(V)
\bigr)^{\top}\bigr] \bigr)^{-1} \bigl(U-\mu_1(V)
\bigr)\elf(\ve)
\\
&&\qquad = \frac{\delta}{E[\delta]} \bigl( J_f E\bigl[\bigl(U-
\mu_1(V)\bigr) \bigl(U-\mu_1(V)\bigr)^{\top} |
\delta=1\bigr] \bigr)^{-1} \bigl(U-\mu_1(V)\bigr)\elf(\ve).
\end{eqnarray*}
This canonical gradient is the influence function of an efficient
estimator of $\vt$. Now use the fact that $\mu_1(V)$ equals
$\mu_{\tilde G}(V)$ and $E[(U-\mu_1(V))(U-\mu_1(V))^{\top}|\delta=1]$
equals $W_{\tilde G}$ to see that this is indeed the characterization
(\ref{eff-mar}).
\end{pf}

%
%re3.2 #&#
\begin{remark}\label{remremextn}
The above efficiency result extends in a straightforward
manner to the case when $V$ is higher dimensional. It also extends to the
partially linear additive model %. In its simplest form, the latter is
%given by
%
\[
Y= \vt^{\top} U + \rho_1(V_1)+
\rho_2(V_2) + \ve,
\]
where $(V_1,V_2)$ takes values in the unit square $[0,1]^2$
and has a density that is bounded and bounded away from zero on
the unit square.
Let $G$ now denote the joint distribution of $(U, V_1, V_2)$.
Assume that the matrix $E[(U-\mu_G(V_1,V_2))(U-\mu_G(V_1,V_2))]^{\top}$,
with $\mu_G(V_1,V_2) = E (U|(V_1,V_2) ),$
is positive definite, and that
$\pi$ is bounded away from zero.
In the present model the space $\mathscr{V}_1$ consists of functions
of the form
\[
a^{\top} U \ell_f(\ve)+ \bigl(b_1(V_1)+b_2(V_2)
\bigr)\ell_f(\ve) + c(\ve),
\]
where $E[b_1^2(V_1)+b_2^2(V_2)]$
is finite. % are square-integrable.
The role of $g$ is now played by
\[
g(X,\delta Y,\delta)= \delta\bigl(U-\tilde\nu_1(V_1)-
\tilde\nu_2(V_2)\bigr) \ell_f(\ve),
\]
where $\tilde\nu_1(V_1)+\tilde\nu_2(V_2)$ minimizes
$ E[\|U-B_1(V_1)-B_2(V_2)\|^2 |\delta=1]$
with respect to functions $B_1$ and $B_2$ from $[0,1]$ into $\R^m$
such that
$E[\|B_1(V_1)\|^2]$ and $E[\|B_2(V_2)\|^2]$ are finite.
The efficient influence function is
\begin{eqnarray*} &&\frac{\delta}{E[\delta]} \bigl(J_f E\bigl[
\bigl(U-\tilde\nu_1(V_1) -\tilde\nu_2(V_2)
\bigr) \bigl(U-\tilde\nu_1(V_1) -\tilde
\nu_2(V_2)\bigr)^{\top}|\delta=1\bigr]
\bigr)^{-1}
\\
&&\qquad{} \times\bigl(U-\tilde\nu_1(V_1)-\tilde
\nu_2(V_2)\bigr)\ell_f(\ve).
\end{eqnarray*}
By the transfer principle, this is the influence function of a complete case
version of an estimator with influence function
\begin{eqnarray*}
&&\bigl(J_f E\bigl[\bigl(U-\nu_1(V_1)-
\nu_2(V_2)\bigr) \bigl(U-\nu_1(V_1)-
\nu_2(V_2)\bigr)^{\top}\bigr]
\bigr)^{-1}\\
&&\qquad{}\times \bigl(U-\nu_1(V_1)-
\nu_2(V_2)\bigr)\ell_f(\ve)
\end{eqnarray*}
in the full model, where
$\nu_1(V_1)+\nu_2(V_2)$ minimizes $E[\|U-B_1(V_1)-B_2(V_2)\|^2]$
over functions $B_1$ and $B_2$ as above.
\citet{Sch96a} constructed estimators in the full model
that have the latter influence function. In particular, he
established their efficiency by showing that the above influence
function is indeed the efficient influence function in the full model.\vadjust{\goodbreak}
\end{remark}
%
%re3.3 #&#
\begin{remark}\label{rem6b}
In the above we have shown that for the partially linear MAR model
(with a possibly additive smooth function) % $\rho$)
an efficient estimator of the parameter $\vt$ can be obtained as the
complete case version of an efficient estimator in the full model.
This is typically also true for other more general semiparametric
regression models and can be verified along the above lines.
%again using the results of M\"uller {\em et al.\ } (2006) and Schick
%(1993).
We sketch this for the following single index model.

In this model $Y=\rho(V+\vt^{\top}U)+\ve$
with one-dimensional $V$, $m$-dimensional $U$, $\vt\in\R^m$ and
twice continuously differentiable $\rho$. Assume again that $\pi$ is
bounded away from zero.
The space $\mathscr{V}$ %$\V$
for this model consists of functions
\[
a^{\top}U \rho'\bigl(V+\vt^{\top}U\bigr)\elf(\ve)
+ b\bigl(V+\vt^{\top}U\bigr) \elf (\ve) + c(\ve)
\]
with $a\in\R^m$, $E[b^2(V+\vt^{\top}U)]<\infty$ and $c$ as before.
%This requires
For this we must require
that $E[\|U\|^2 (\rho'(V+\vt^{\top}U))^2]$ is finite.
Now one works with $g(X,\delta Y, \delta)= \delta(U-\nu_1(V+\vt^{\top}U))
\rho'(V+\vt^{\top}U) \elf(\ve)$ and $\nu_1(V+\vt^{\top}U)=
E(U|V+\vt^{\top} U,\delta=1)$, and
%verifies that the canonical gradient is
obtains the canonical gradient
\[
g_*(X,\delta Y,\delta)= \frac{\delta}{E[\delta]} (J_f W_1)^{-1}
\bigl(U-\nu_1\bigl(V+\vt^{\top}U\bigr)\bigr)
\rho'\bigl(V+\vt^{\top}U\bigr) \elf (\ve)
\]
if $W_1=E[(U-\nu_1(V+\vt^{\top}U))(U-\nu_1(V+\vt^{\top}U))^{\top
} (\rho'(V+
\vt^{\top}U))^2|\delta=1]$ is invertible.
By the transfer principle, this is
the influence function of a complete case version
of an estimator with influence function
$(J_fW)^{-1} (U-\nu(V+\vt^{\top}U))\rho'(V+\vt^{\top}U)\elf(\ve)$,
where $\nu(V+\vt^{\top}U)= E[U|V+\vt^{\top}V]$ and
$W=E[(U-\nu(V+\vt^{\top}U))(U-\nu(V+\vt^{\top}U))^{\top}](\rho'(V+\vt^{\top}U))^2]$.
The latter influence function is the efficient gradient for the full model.
Indeed, it is the canonical gradient for the case when $\delta=1$
almost surely.
\end{remark}

%s4 #&#
\section{Testing for normal errors}
\label{sec4}

In this section we shall introduce a test for normal errors
which uses the Khmaladze transform of the
empirical distribution function
$\hemp(t)= n^{-1} \sum_{j=1}^n\1[\hve_j \leq t]$, $t \in\R$,
based on residuals $\hve_j$.
Goodness-of-fit tests for the full model based on that transform
were discussed in Khmaladze and Koul (\citeyear{khkou04}, \citeyear
{khkou09}) for parametric
and nonparametric regression, and by MSW for the partially
linear regression model considered here. Due to the transfer principle,
it is now straightforward to adapt the approach by MSW to the MAR model,
which is what we will do here for a simple illustration of the method.
Note that MSW consider the more complex case where $V$ is a covariate
\textit{vector}.

First, we briefly sketch the approach for the full model.
To avoid additional assumptions, we estimate $\vt$ and $\rho$ using
a least squares approach with the trigonometric basis. This is discussed
in MSW, Section 4, for an additive regression
function, that is, with $\rho(x_1,\ndots,x_q)=\rho_1(x_1) + \cdots+
\rho_q(x_q)$.
(Here we have $q=1$.) For $k=1,2,\ndots,$ we set
\[
\phi_k(x)= \cos(\pi k x),\qquad  0\le x \le1.\vadjust{\goodbreak}
\]
Our estimator of the regression function $r(u,v)=\vt^{\top}u + \rho(v)$
is then
\[
\hat r(u,v)= \hvt^{\top}u + \sum_{k=0}^K
\hat\beta_k \phi_k(v),
\]
where $\phi_0(x)=1$ and $(\hvt^{\top},\hat\beta_0,\ndots, \hat
\beta_K)$
minimizes
\[
\sum_{j=1}^n \Biggl( Y_j -a^{\top} U_j - \sum_{k=0}^K b_k \phi_k(V_j) \Biggr)^2
\]
with respect to $a, b_0, \ndots, b_K$.
For $j=1,\ndots,n$ the error $\ve_j$ is estimated by the residual
\[
\hve_j = Y_j - \hvt^\top U_j -
\sum_{k=0}^K \hat\beta_k
\phi_k(V_j).
\]
We also need the normalized residuals $\hat Z_j= \hve_j/\hs$,
where $\hs$ is the square root of $(1/n)\sum_{j=1}^n\hve_j^2$.

Assume for the remainder of this section that $f$ has finite
Fisher information for location and finite fourth moment.
This assumption is met by the normal density.
It then follows from MSW, Theorem 4.1 and Remark 4.2,
that with $K=K_n \sim n^{-1/4}$ we have the uniform stochastic expansions
%
%e4.1 #&#
\begin{equation}
\label{uu1} \sup_{t\in\R} \Biggl| \frac{1}{\sqrt{n}} \sum
_{j=1}^n \bigl(\1[\hve_j \leq t] - \1[
\ve_j \leq t] - f(t) \ve_j \bigr) \Biggr| =
o_P(1)
\end{equation}
and
%
%e4.2 #&#
\begin{equation}
\label{uu2} \qquad \sup_{t\in\R} \Biggl|\frac{1}{\sqrt{n}} \sum
_{j=1}^n \biggl(\1[\hat Z_j\le t] -
\1[Z_j\le t] - f_*(t) \biggl(Z_j + t\frac{Z_j^2-1}{2}
\biggr) \biggr) \Biggr| = o_P(1),
\end{equation}
where $f_*$ denotes the density of the normalized errors
$Z_j = \ve_j/\sigma$.

Write $\phi$ for the standard normal density.
In terms of the density $f_*$, the null hypothesis is
\[
H_0\dvtx f_*=\phi.
\]
MSW proposed the test statistic
\[
T_{n}= \sup_{t\in\R}\Biggl | \frac{1}{\sqrt{n}} \sum
_{j=1}^n \bigl(\1[\hat Z_j \le t] - H(t
\wedge\hat Z_j) h(\hat Z_j) \bigr) \Biggr|,
\]
with
%
%e4.3 #&#
\begin{eqnarray}
\label{hH} h(x) & =& \bigl(1, x, x^2-1\bigr)^{\top},\qquad
\Gamma(x) = \int_{x}^{\infty} h(z)h^{\top}(z)
\phi(z) \,dz,
\nonumber
\\[-8pt]
\\[-8pt]
\nonumber
H(t) & =& \int_{-\infty}^t h^{\top}(x)
\Gamma^{-1}(x) \phi(x) \,dx.
\end{eqnarray}
This is a version of the martingale transform test
of \citet{khkou09}
for fitting an error distribution in nonparametric regression.
MSW showed that under the null hypothesis the test statistic $T_n$
converges in distribution to $\zeta$, which is the supremum of a
standard Brownian motion given in (\ref{z}).
This holds for every distribution function $G$ satisfying (G1) and (G2).
Since $\ve$ and $(\delta,U,V)$ are independent under the MAR assumption,
the conditional distribution of $(\ve,U,V)$, given $\delta=1$,
is given by $F\times\tilde G$, where $\tilde G$ is the conditional distribution
of $(U,V)$, given $\delta=1$. Thus, if $\tilde G$ satisfies (G1) and (G2),
then the transfer principle applies and yields the same limiting
distribution for the complete case version $T_{c}$ of $T_n$, where
%
%e4.4 #&#
\begin{eqnarray}
\label{tc} T_{c}= \sup_{t\in\R} \Biggl| \frac{1}{\sqrt{N}} \sum
_{j=1}^n \delta_j \bigl(\1[
\hat Z_{jc} \leq t] - H(t\wedge\hat Z_{jc}) h(\hat
Z_{jc}) \bigr) \Biggr|.
\end{eqnarray}
Here $\hat Z_{jc}$ are the complete case versions
of the normalized residuals
and are defined by $\hve_{jc}/\hat\sigma_{c}$ with
$\hve_{jc} = Y_j - \hvt_c^\top U_j - \sum_{k=0}^{K_N} \hat\beta_k
\psi_k(V_j)$ and
$\hat\sigma_c$ the square root of $N^{-1} \sum_{j=1}^n\dj\hve_{jc}^2$,
while $(\hvt_c^{\top},\hat\beta_0,\ndots,\hat\beta_{K_N})$
are the least squares estimators minimizing
\[
\sum_{j=1}^n\dj \Biggl( Y_j
- a^{\top} U_j - \sum_{k=0}^{K_N}
b_k \phi_k(V_j) \Biggr)^2.
\]

The transfer principle for asymptotically linear statistics also provides
complete case versions of the expansions (\ref{uu1}) and (\ref{uu2})
from above.
The first expansion becomes
\[
\sup_{t\in\R} \Biggl| \frac{1}{\sqrt{N}} \sum_{j=1}^n
\delta_j \bigl(\1[\hve_{jc} \leq t] - \1[
\ve_j \leq t] - f(t) \ve_j \bigr) \Biggl| =
o_P(1),
\]
and the second expansion becomes\vspace*{-1pt}
\[
\sup_{t\in\R} \Biggl| \frac{1}{\sqrt{N}} \sum_{j=1}^n
\delta_j \biggl(\1[\hat Z_{jc} \le t] -
\1[Z_j\le t] - f_*(t) \biggl(Z_j + t\frac{Z_j^2-1}{2}
\biggr) \biggr)\Biggr | = o_P(1).\vspace*{-2pt}
\]

%s5 #&#
\section{Testing for linearity}
\label{sec5}

In this section we address testing whether the function $\rho$ in the
partially linear MAR model is constant. In the previous
section we demonstrated how the transfer principle can be used
to adapt a known test for the full model to the MAR model. We now
show how to develop a test procedure for the MAR model
when no counterpart to the full model exists. Our approach is to first
develop a procedure for the full model, and then to apply the
transfer principle. Our test statistic is inspired by that in
Stute, Xu and Zhu (\citeyear{stu08}).

Under the null hypothesis the partially linear model reduces
to the linear regression model $Y= \alpha+\vt^{\top}U +\ve$, where
$\alpha$ is an unknown constant,
that is, we have\vspace*{-1pt}
\[
H_0\dvtx \rho(v) = \alpha\qquad\mbox{for all $v \in\R$ and some $\alpha
\in\R$.}\vspace*{-1pt}
\]
To simplify notation, we introduce $\beta=(\alpha,\vt^{\top})^{\top}$
and $Z=(1,U^{\top})^{\top}$. Then we can write the model under
the null hypothesis as $Y=\beta^{\top}Z+\ve$.

It follows from (G2) that the dispersion matrix $\Lambda_G$ of $U$
is positive definite. From this we immediately see that the matrix\vspace*{-1pt}
\[
M_G= E\bigl[ZZ^{\top}\bigr] = \left[\matrix{ 1 & E
\bigl[U^{\top}\bigr] \vspace*{2pt}
\cr
E[U] & E\bigl[UU^{\top}
\bigr]} \right]\vspace*{-1pt}
\]
is also positive definite. Thus, the least squares estimator
$\hat\beta$ of $\beta=(\alpha,\vt^{\top})^{\top}$
is root-$n$ consistent under the null hypothesis, as it satisfies\vspace*{-1pt}
\[
\hat\beta= \beta+ M_G^{-1} \frac1n \sum
_{j=1}^n Z_j \ve_j +
o_P\bigl(n^{-1/2}\bigr).\vspace*{-1pt}
\]

Now let $\h$ denote a continuous non-constant function on $[0,1]$.
Introduce the least squares estimator $\hat\gamma$ for regressing
the responses $\h(V_j)$ on the design vectors~$Z_j$,
so that $\hat\gamma$ minimizes\vspace*{-1pt}
\[
\frac1n \sum_{j=1}^n \bigl(
\h(V_j)-\gamma^{\top}Z_j\bigr)^2.\vspace*{-1pt}
\]
Set $R_j = \h(V_j)-\hat\gamma^{\top} Z_j$, $W_j= R_j/(n^{-1} \sum_{j=1}^n
R_j^2)^{1/2}$, and
$\hve_{j0}= Y_j-\hat\beta^{\top}Z_j$, $j=1,\ndots,n$.
Our test statistic in the full model is\vspace*{-1pt}
\[
T_n = \sup_{t\in\R} \Biggl|\frac{1}{\sqrt{n}} \sum
_{j=1}^n W_j \1 [\hve_{j0}\le
t] \Biggr|.\vspace*{-1pt}
\]
As in Stute, Xu and Zhu (\citeyear{stu08}), we have the following result.

%le5.1 #&#
\begin{lemma}\label{lm3} Suppose the null hypothesis holds and $f$ is
uniformly continuous.
Then $T_n$ converges in distribution to
$\zeta_0 = \sup_{0\le t \le1} |B_0(t)|$,
where $B_0$ denotes a standard Brownian bridge.
\end{lemma}

\begin{pf}
Set
\[
\h_G(X)=\h(V)-\gamma_G^{\top}Z,
\]
where $\gamma_G$ minimizes $E[(\h(V)-\gamma^{\top}Z)^2]$.
Let $\rho_G=\Lambda_G^{-1}E[\h(V)(U-E[U])]$. Then it is easy to
check that
\begin{eqnarray*} \h_G(X) &=&\h(V)-E\bigl[\h(V)\bigr]-
\rho_G^{\top}\bigl(U-E[U]\bigr)
\\
&=&\h(V)-E\bigl[\h(V)\bigr]-\rho_G^{\top} \bigl(
\mu_G(V)-E[U]\bigr) - \rho_G^{\top
} \bigl(U-
\mu_G(V)\bigr).
\end{eqnarray*}
Note also that $\h$ being non-constant on $[0,1]$ and $V$ having a
positive density on $[0,1]$ implies
$\h(V)$ has a positive variance. These facts together with $W_G$ being
positive definite guarantee that
$E[\h^2_G(X)]=\var(\h(V)-\rho_G^{\top} \mu_G(V))+ \rho_G^{\top}
W_G\rho_G$
is positive.

Next, let $g$ be a measurable function such that $E[g^2(X)]$ is finite and
assume $f$ is uniformly continuous. Then Theorem 2.2.4 of \citet
{kou02} yields
\[
\sup_{t\in\R} \Biggl| \frac1n \sum_{j=1}^n
g(X_j) \bigl(\1[\hve_{j0} \le t] - \1[\ve_j
\le t]\bigr) -f(t) E\bigl[g(X)Z^{\top}\bigr](\hat\beta-\beta)\Biggr | =
o_P\bigl(n^{-1/2}\bigr).
\]
From this fact we obtain
\[
\sup_{t\in\R} \Biggl| \frac1n \sum_{j=1}^n
R_j\bigl(\1[\hve_{j0} \le t]-\1[\ve_j\le t]
\bigr) -f(t) \hat D(\hat\beta-\beta) \Biggr|= o_P\bigl(n^{-1/2}
\bigr),
\]
where
\[
\hat D= E\bigl[\h(V)Z^{\top}\bigr]- \hat\gamma^{\top} E
\bigl[ZZ^{\top}\bigr] = E\bigl[\h_G(V)Z^{\top}\bigr]+
o_P(1).
\]
In view of the identities $E[\h_G(V)Z^{\top}]=0$ and $\sum_{j=1}^nR_j=0$,
we can conclude
\[
\sup_{t\in\R} \Biggl| \frac1n \sum_{j=1}^n
R_j\1[\hve_{j0} \le t]-\frac1n \sum
_{j=1}^n R_j\bigl(\1[\ve_j
\le t]-F(t)\bigr) \Biggr|=o_P\bigl(n^{-1/2}\bigr).
\]
Writing $R_j-\h_G(V_j)= - (\hat\gamma-\gamma_G)^{\top} Z_j$,
we derive the expansions
\begin{eqnarray*}
\sup_{t\in\R} \Biggl|\frac1n \sum_{j=1}^n
\bigl(R_j-\h_G(V_j)\bigr) \bigl(\1[
\ve_j\le t]-F(t)\bigr) \Biggr| &=&o_P\bigl(n^{-1/2}
\bigr),
\\
\frac1n \sum_{j=1}^n
\bigl(R_j-\h_g(V_j)\bigr)^2 \le
\|\hat\gamma-\gamma_G\|^2 \frac1n \sum
_{j=1}^n \|Z_j\|^2 &=&
o_P(1),
\end{eqnarray*}
and therefore obtain $n^{-1} \sum_{j=1}^nR_j^2 = E[\h_G^2(V)] + o_P(1)$.
The above derivations in turn yield
\[
\sup_{t\in\R} \Biggl| \frac{1}{\sqrt{n}} \sum_{j=1}^n
W_j\1[\hve_{j0}\le t] - \frac{1}{\sqrt{n}} \sum
_{j=1}^n \h^*_G(V_j)
\bigl(\1[\ve_j\le t]-F(t)\bigr) \Biggr| = o_P(1),
\]
with $\h_G^*= \h_G/E[\h^2_G(V)]^{1/2}$.
Since, again by Theorem 2.2.4 of \citet{kou02}, the process
\[
\frac{1}{\sqrt{n}} \sum_{j=1}^n
\h^*_G(V_j) \bigl(\1[\ve_j\le t]-F(t)
\bigr),\qquad -\infty\le t\le \infty,
\]
converges in $D([-\infty,\infty])$ to a time-changed Brownian bridge $B_0(F)$,
we conclude that $T_n$ has the desired limiting distribution.
\end{pf}

The complete case version of $T_n$ is given by
%
%e5.1 #&#
\begin{equation}
\label{tc2} T_{c} = \sup_{t\in\R} \Biggl|\frac{1}{\sqrt{N}} \sum
_{j=1}^n \delta_jR_{jc}
\1[\he_{jc} \le t] \Biggr| \bigg/ \Biggl(\frac1N \sum
_{j=1}^n\dj R_{jc}^2
\Biggr)^{1/2},
\end{equation}
with $\he_{jc} = Y_j-\hat\beta_c^{\top} Z_j$,
$R_{jc}=\h(V_j)-\hat\gamma_{c}^{\top} Z_j$,
\begin{eqnarray*}
\hat\beta_c&=& \argmin_{b} \sum
_{j=1}^n\dj\bigl(Y_j-b^{\top}
Z_j\bigr)^2\quad \mbox{and}\\
 \hat\gamma_{c}&=&
\argmin_{\gamma} \sum_{j=1}^n\dj
\bigl(\h(V_j)-\gamma^{\top} Z_j
\bigr)^2 .
\end{eqnarray*}
By the transfer principle, the limiting distribution of $T_{c}$
under the null hypothesis will be that of
$\zeta_0 $ %of (\ref{z0}),
from the above lemma, as long as $f$ is uniformly continuous and~$\tilde G$ satisfies (G1) and (G2).

%re5.1 #&#
\begin{remark}\label{remgen}
The above is easily extended to cover testing for other parametric forms
for $\rho$. For example, we can test whether $\rho$ is linear,
$\rho(v)=a+bv$. In this case we proceed as above, but with the role of $Z$
now played by the vector $(1,U^{\top},V)$ and with $\chi$ chosen to be
nonlinear.
\end{remark}

\section*{Acknowledgments}
The authors would like to
thank the two reviewers for their constructive comments,
which have helped to improve the paper.

% imsref loaded by akundreckaite, 2012-11-22 13:52:18
%

%suskaldyti doi

\printaddresses

\end{document}